\documentclass[11pt]{article}
\usepackage{amsmath}
\usepackage{amssymb}
\usepackage{amsthm}
\usepackage{comment}
\usepackage[morefloats=200]{morefloats}
\usepackage{enumitem}
\setlist{nosep}

\newcommand{\Z}{{\mathbb Z}}
\newtheorem{theorem}{Theorem}[section]

\newtheorem{lemma}[theorem]{Lemma}

\newtheorem{example}[theorem]{Example}

\begin{document}
\baselineskip 18pt
\title{Direct constructions of some group divisible designs with block~size~$4$ and up to $50$~points}

{\small
\author
{R. Julian R. Abel \\
 School of Mathematics and Statistics\\
  UNSW Sydney\\ NSW 2052, Australia\\
 \texttt{r.j.abel@unsw.edu.au }
 \and
   Thomas Britz  \\
  School of Mathematics and Statistics\\
   UNSW Sydney\\ NSW 2052, Australia\\
  \texttt{britz@unsw.edu.au  }
 \and
   Yudhistira A.  Bunjamin  \\
  School of Mathematics and Statistics\\
   UNSW Sydney\\ NSW 2052, Australia\\
  \texttt{yudhi@unsw.edu.au  }
  \and
  Diana Combe  \\
  School of Mathematics and Statistics\\
   UNSW Sydney\\ NSW 2052, Australia\\
  \texttt{diana@unsw.edu.au }
}    % end author
}    % end small

\date{}

\maketitle 

\noindent{\bf Abstract:} 
In this note, we give direct constructions of some group divisible designs (GDDs) with block size $4$ that have up to $50$ points.

\noindent{\bf Keywords:} group divisible design (GDD), feasible group type.

\noindent{\bf Mathematics Subject Classification:}  05B05

\section{Introduction}

Let $X$ be a finite set of {\it points} with a partition into parts which we call {\it groups}. 
Any $k$-element subset of $X$ is called a {\it block}. 
A collection of {\it blocks} is a {\it group divisible design} with block size $k$, or $k$-GDD if 
(i) no two points from the same group appear together in any block and 
(ii) any two points from distinct groups appear together in exactly one block. 
The {\it group type} (or {\it type}) of a $k$-GDD is the multiset $\left\{|G|: G {\text{ is a part of the partition}}\right\}$. 
The group type can also be expressed in  `exponential' notation where the type $t_1^{u_1}  t_2^{u_2} \ldots t_m^{u_m}$ means there are $u_i$ groups of size $t_i$ for $i=1,2, \dots, m$.

There are known necessary conditions for the existence of a $4$-GDD of type $\{g_1, g_2, \ldots, g_m\}$.
These are given in \cite{ABC.50less}.  
However, these necessary conditions are not sufficient. 

In this note, we are concerned with $4$-GDDs. 
We say that a multiset $\{g_1, g_2, \ldots, g_m\}$ of positive integers is a {\it feasible} group type for a $4$-GDD if it satisfies the known necessary conditions. 

In~\cite{krestin}, Kreher and Stinson gave the list of feasible group types for $4$-GDDs with up to $v=30$ points. 
There the existence question was answered for all but three group types.
The existence results for those three types were completed in~\cite{ABBC.2t5s, ABC.30less}.

In \cite{ABC.50less}, the list of feasible group types for $4$-GDDs was extended to those with up to $v=50$ points.
For $31 \leq v \leq 50$ and $v \equiv 0 \pmod{3}$ there exist $4$-GDDs for all feasible group types; this result was completed in~\cite{ABC.50less}. 
In~\cite{ABBCF.1mod3}, the feasible types for $31 \leq v \leq 50$ and $v \equiv 1$ and $2 \pmod{3}$ were considered.  
The questions of existence were completed for $v \equiv 1 \pmod{3}$; and for $v \equiv 2 \pmod{3}$ the known results were extended leaving unknown the question of existence for types $2^{11} 8^1 11^1$, $2^1 5^4 8^1 11^1$, $2^6 5^2 11^2$,
$2^{5} 5^3 8^1 11^1$, 
$2^{2} 5^2 8^1 11^2$,  $2^{1} 5^3 8^2 11^1$, 
 $2^{5} 5^3 11^2$, $2^{2} 5^2 11^3$, $2^{1} 5^3 8^1 11^2$, $5^{4} 8^2 11^1$,
  $2^{9} 5^2 11^2$, 
 $2^{8} 5^3 8^1 11^1$,   $2^6 5^1 11^3$, $2^{5} 5^2 8^1 11^2$,  $2^{4} 5^4 8^1 14^1$, $2^{4} 5^3 8^2 11^1$,  $2^{3} 11^4$,  
$2^{2} 5^7 11^1$,  $2^{2} 5^1 8^1 11^3$,  $2^{1} 5^2 8^2 11^2$ and   $5^{3} 8^3 11^1$.

In this note, we provide direct constructions of $4$-GDDs with up to $50$ points, particularly for the group types for which the question of existence was left open in~\cite{ABBCF.1mod3}.
The primary intention of this note is to make these constructions publicly available before they appear in future publications by the current authors.
Occasionally, we may also include constructions of $4$-GDDs of group types for which existence has already been determined.
However, these constructions will be different from the ones that were first used to prove the existence of $4$-GDDs of these group types.
As a result, this note may be updated from time to time with new constructions.

\section{Direct constructions by assuming an automorphism}
\label{s:direct.constructions}

The $4$-GDDs in this section have been found directly by assuming the existence of a cyclic automorphism group, $\mathbb{G}$. 
%See  for instance~\cite{ABBC.2t8s,ABBC.2t5s,ABBCF.1mod3,ABC.30less,ABC.50less,ABC.3562,Forbes2,Forbes3,Forbesgcc,Forbes1,WeiGe.4gdd.gpn1.g0mod6,gegum}.  
For each $4$-GDD, a set of base blocks is given.
%For the $4$-GDDs constructed directly in this paper, the point set of the design consists of one or more copies of the group~$\mathbb{G}$ (here, $\mathbb{G}$ is $\mathbb{Z}_5$, $\mathbb{Z}_6$, $\mathbb{Z}_7$ or $\mathbb{Z}_{14}$) plus possibly a few copies of $\mathbb{Z}_2$ or $\mathbb{Z}_3$ and possibly one or more infinite points.  
The blocks are obtained by developing the subscripts of the points in each base block from copies of $\mathbb{G}$ over $\mathbb{G}$;
the infinite points remain unaltered when developed.  
%When the point set includes any extra copies of $\mathbb{Z}_2$ or $\mathbb{Z}_3$, those points are developed over $\mathbb{Z}_2$ or $\mathbb{Z}_3$ as the others are developed over~$\mathbb{G}$.
%Also there are usually a few blocks which remain invariant when some nonzero element of $\mathbb{G}$ is added to it;
%the number of blocks generated by any one of those blocks is less than the size of $\mathbb{G}$.

%The $4$-GDDs below appear in the order in which they were added to this note.

We note that a $4$-GDD of type~$3^8 6^1 12^1$ was already known to exist prior to this note.
The $4$-GDD of type~$3^8 6^1 12^1$ in \cite{reesstinsubs} was obtained by Rees and Stinson implicitly.
A $\{3,4\}$-GDD is defined similarly to a $k$-GDD except that each block is either of size~$3$ or size~$4$.
Rees and Stinson construct a $\{3,4\}$-GDD of type~$3^8 6^1$ whose triples fall into $12$~resolution classes in \cite[Appendix]{reesstinsubs}. 
To construct a $4$-GDD of type~$3^8 6^1 12^1$, they started with the $\{3,4\}$-GDD of type~$3^8 6^1$ and turned the blocks of size~$3$ into blocks of size~$4$ by adding each point from the group of size~$12$ to the blocks of size~$3$ in one of the resolution classes.

In Example~\ref{example.4gdd.3^8_6^1_12^1}, we give a $4$-GDD of type~$3^8 6^1 12^1$ that was constructed directly.

\begin{example} \rm
\label{example.4gdd.3^8_6^1_12^1}
A $4$-GDD of type~$3^8 6^1 12^1$ can be constructed in the following manner.

The points in the $4$-GDD are
\begin{itemize}
    \item $a_i,b_i,q_i,r_i,y_i,z_i$ for $i\in\Z_6;$
    \item $p_i$ for $i\in\Z_3;$ and
    \item $\infty_1,\infty_2,\infty_3$.
\end{itemize}
The groups are
\begin{itemize}
    \item $\{a_i: i \in \Z_6\} \cup  \{b_i: i \in \Z_6\};$
    \item $\{p_i: i \in \Z_3\} \cup \{ \infty_1, \infty_2,\infty_3\};$
    \item $\{q_i, q_{i+2}, q_{i+4}\}$, $\{r_i, r_{i+2}, r_{i+4}\}$, $\{y_i, y_{i+2}, y_{i+4}\};$ and
    \item $\{z_i, z_{i+2}, z_{i+4}\}$ for $i \in \{0,1\}$.
\end{itemize}
Blocks are obtained by developing modulo $6$ the subscripts of the noninfinite points in the base blocks given in Table~\ref{table.4gdd.3^8_6^1_12^1}.
The first two blocks in the first column generate three blocks each.
\end{example}

\begin{table}[ht]
    \centering
    \caption{Base blocks for the $4$-GDD of type~$3^8 6^1 12^1$ in Example~\ref{example.4gdd.3^8_6^1_12^1}}
    \label{table.4gdd.3^8_6^1_12^1}
    \renewcommand{\arraystretch}{1.2}
    \[ \begin{array}{|l|l|l|l|l|l|}\hline
    \{p_0,p_3,q_0,q_3\} & \{a_0,p_0,p_1,y_0\}      &  \{a_0,q_3,u_5,\infty_3\}      &    \{b_0,p_4,q_5,u_0\}          &  \{b_0,p_5,u_4,u_5\}    \\\hline
    \{t_0,t_3,u_0,u_3\} & \{a_0,p_2,q_0,t_0\}      &  \{a_0,q_1,t_3,u_1\}           &   \{b_0,p_0,q_2,\infty_2\}      &  \{b_0,q_0,q_1,t_4\}    \\\hline
                        & \{a_0,p_3,q_2,u_0\}      &  \{a_0,q_4,t_5,y_2\}           &   \{b_0,p_2,u_3,y_0\}           &  \{b_0,q_4,t_3,y_1\}     \\\hline
                        & \{a_0,p_4,t_1,\infty_1\} &  \{a_0,q_5,u_2,y_1\}           &  \{b_0,p_1,t_1,t_2\}            &  \{b_0,q_3,u_2,\infty_1\}  \\\hline
                        & \{a_0,p_5,t_4,u_3\}      &  \{a_0,t_2,u_4,\infty_2\}      &  \{b_0,p_3,t_5,\infty_3\}       &  \{b_0,t_0,u_1,y_2\}        \\\hline
  \end{array} \]
\end{table}

\begin{lemma}
\label{lemma.4gdd.2^6_5^2_11^2}
There exists a $4$-GDD of type~$2^6 5^2 11^2$.

\begin{proof}
The points in the $4$-GDD are
\begin{itemize}
    \item $a_i,b_i,c_i,f_i,g_i,p_i,q_i,r_i,s_i,t_i$ for $i\in\Z_4$; and
    \item $h_i,u_i$ for $i\in\Z_2$.
\end{itemize}
The groups are
\begin{itemize}
    \item $\{a_i,b_i\}$ for $i\in\Z_4$;
    \item $\{c_i,c_{i+2}\}$;
    \item $\{f_i,f_{i+2},g_i,g_{i+2},h_i\}$; and
    \item $\{p_i,p_{i+2},q_i,q_{i+2},r_i,r_{i+2},s_i,s_{i+2},t_i,t_{i+2},u_i\}$ for $i\in \{0,1\}$.
\end{itemize}
Blocks are obtained by developing modulo $4$ the subscripts of the noninfinite points in the base blocks given in Table~\ref{table.4gdd.2^6_5^2_11^2}.
The first block in the first column generates one block, while the second block in that column generates two blocks.
\end{proof}
\end{lemma}

\begin{table}[ht]
    \centering
    \caption{Base blocks for the $4$-GDD of type~$2^6 5^2 11^2$ in Lemma~\ref{lemma.4gdd.2^6_5^2_11^2}}
    \label{table.4gdd.2^6_5^2_11^2}
    \renewcommand{\arraystretch}{1.2}
    \[ \begin{array}{|l|l|l|l|l|l|}\hline
    \{h_{0},h_{1},u_{0},u_{1}\} & \{a_{0},a_{1},p_{0},r_{1}\} & \{a_{0},g_{2},r_{3},s_{2}\} &  \{b_{0},f_{3},t_{0},u_{1}\} & \{c_{0},g_{1},r_{3},u_{0}\} \\\hline
    \{a_{0},a_{2},b_{1},b_{3}\} & \{a_{0},b_{2},c_{0},h_{0}\} & \{a_{0},g_{3},r_{2},t_{1}\} &  \{b_{0},g_{2},p_{0},s_{3}\} & \{c_{0},h_{1},r_{2},s_{3}\} \\\hline
                                & \{a_{0},c_{1},f_{1},g_{0}\} & \{a_{0},h_{1},q_{1},t_{0}\} &  \{b_{0},g_{3},p_{3},t_{2}\} & \{f_{0},f_{1},p_{3},s_{0}\} \\\hline
                                & \{a_{0},c_{2},p_{1},p_{2}\} & \{b_{0},b_{1},q_{2},s_{1}\} &  \{b_{0},g_{0},r_{0},t_{1}\} & \{f_{0},h_{1},p_{0},r_{3}\} \\\hline
                                & \{a_{0},c_{3},s_{0},s_{3}\} & \{b_{0},c_{1},f_{0},g_{1}\} &  \{b_{0},h_{1},s_{2},t_{3}\} & \{g_{0},g_{1},q_{2},s_{3}\}  \\\hline
                                & \{a_{0},f_{2},q_{2},q_{3}\} & \{b_{0},c_{3},p_{1},u_{0}\} &  \{c_{0},c_{1},q_{3},t_{0}\} & \{g_{0},h_{1},p_{1},q_{0}\} \\\hline
                                & \{a_{0},f_{3},s_{1},u_{0}\} & \{b_{0},c_{0},q_{0},r_{1}\} &  \{c_{0},f_{2},q_{1},r_{0}\} &   \\ \hline
                                & \{a_{0},f_{0},t_{2},t_{3}\} & \{b_{0},f_{1},p_{2},q_{3}\} &  \{c_{0},f_{1},s_{2},t_{1}\} &   \\ \hline
                                & \{a_{0},g_{1},q_{0},u_{1}\} & \{b_{0},f_{2},r_{2},r_{3}\} &  \{c_{0},g_{2},p_{1},t_{2}\} &  \\ \hline 
 \end{array} \]
\end{table}

\begin{lemma}
\label{lemma.4gdd.2^9_5^2_11^2}
There exists a $4$-GDD of type~$2^9 5^2 11^2$.

\begin{proof}
The points in the $4$-GDD are
\begin{itemize}
    \item $a_i,b_i,c_i,d_i,f_i,g_i,p_i,q_i,r_i,s_i,t_i$ for $i\in\Z_4$; 
    \item $h_i,u_i$ for $i\in\Z_2$; and
    \item $\infty_1,\infty_2$.
\end{itemize}
The groups are
\begin{itemize}
    \item $\{a_i,b_i\}$ and $\{c_i,d_i\}$ for $i\in\Z_4$;
    \item $\{f_i,f_{i+2},g_i,g_{i+2},h_i\}$;
    \item $\{p_i,p_{i+2},q_i,q_{i+2},r_i,r_{i+2},s_i,s_{i+2},t_i,t_{i+2},u_i\}$ for $i\in \{0,1\}$; and
    \item $\{\infty_1,\infty_2\}$.
\end{itemize}
Blocks are obtained by developing modulo $4$ the subscripts of the noninfinite points in the base blocks given in Table~\ref{table.4gdd.2^9_5^2_11^2}.
The first block in the first column generates one block, while the other four blocks in that column generate two blocks each.
\end{proof}
\end{lemma}

\begin{table}[ht]
    \centering
    \caption{Base blocks for the $4$-GDD of type~$2^9 5^2 11^2$ in Lemma~\ref{lemma.4gdd.2^9_5^2_11^2}}
    \label{table.4gdd.2^9_5^2_11^2}
    \renewcommand{\arraystretch}{1.2}
    \[ \begin{array}{|l|l|l|l|l|}\hline
    \{h_{0},h_{1},u_{0},u_{1}\}    & \{a_{0},a_{1},c_{0},g_{0}\} & \{a_{0},g_{1},p_{2},r_{1}\}     & \{b_{0},g_{2},g_{3},q_{2}\}    &  \{d_{0},f_{1},p_{0},u_{1}\}  \\\hline
    \{a_{0},a_{2},u_{0},\infty_1\} & \{a_{0},b_{1},f_{0},h_{1}\} & \{a_{0},h_{0},q_{2},s_{3}\}     & \{b_{0},g_{1},t_{1},u_{0}\}    &  \{d_{0},f_{2},r_{1},\infty_2\} \\\hline
    \{b_{0},b_{2},u_{1},\infty_2\} & \{a_{0},b_{2},p_{1},s_{0}\} & \{a_{0},q_{3},s_{2},\infty_2\}  & \{b_{0},h_{1},p_{2},r_{3}\}    &   \{d_{0},g_{2},q_{3},u_{0}\}  \\\hline
    \{c_{0},c_{2},h_{0},\infty_2\} & \{a_{0},b_{3},q_{0},t_{1}\} & \{b_{0},b_{1},c_{0},p_{1}\}     &  \{c_{0},c_{1},f_{0},q_{0}\}   &  \{d_{0},g_{0},s_{1},s_{2}\}  \\\hline
    \{d_{0},d_{2},h_{1},\infty_1\} & \{a_{0},c_{1},r_{0},s_{1}\} & \{b_{0},c_{1},d_{0},t_{0}\}     & \{c_{0},d_{2},g_{1},r_{2}\}    &  \{d_{0},h_{0},p_{2},q_{1}\}  \\\hline
                                   & \{a_{0},c_{2},r_{2},u_{1}\} & \{b_{0},c_{2},g_{0},s_{0}\}     & \{c_{0},d_{1},q_{1},t_{0}\}    &  \{f_{0},f_{1},s_{2},t_{3}\} \\\hline
                                   & \{a_{0},d_{1},f_{1},g_{2}\} & \{b_{0},d_{1},d_{2},s_{1}\}     & \{c_{0},f_{2},p_{0},p_{3}\}    &  \{f_{0},g_{3},q_{1},r_{2}\} \\\hline
                                   & \{a_{0},d_{0},f_{3},p_{3}\} & \{b_{0},d_{3},r_{1},r_{2}\}     & \{c_{0},f_{1},s_{1},u_{0}\}    &  \{g_{0},h_{1},r_{2},t_{3}\}  \\\hline
                                   & \{a_{0},d_{3},p_{0},q_{1}\} & \{b_{0},f_{1},q_{0},q_{3}\}     & \{c_{0},g_{3},p_{2},t_{1}      &  \{g_{0},p_{2},s_{3},\infty_1\} \\\hline
                                   & \{a_{0},d_{2},t_{0},t_{3}\} & \{b_{0},f_{0},r_{0},s_{3}\}     & \{c_{0},h_{1},s_{3},t_{2}\}    &  \{g_{0},p_{0},t_{1},\infty_2\}  \\\hline
                                   & \{a_{0},f_{2},r_{3},t_{2}\} & \{b_{0},f_{2},t_{3},\infty_1\}\ & \{c_{0},q_{2},r_{1},\infty_1\} &     \\\hline

    \end{array} \]
\end{table}

\begin{lemma}
\label{lemma.4gdd.2^6_5^1_11^3}
There exists a $4$-GDD of type~$2^6 5^1 11^3$.

\begin{proof}
The points in the $4$-GDD are
\begin{itemize}
    \item $a_i,b_i,c_i,d_i,e_i,s_i,t_i$ for $i\in\Z_6$;
    \item $f_i$ for $i \in\Z_3$;
    \item $p_i,q_i$ for $i\in\Z_2$; and
    \item $\infty$.
\end{itemize}
The groups are
\begin{itemize}
    \item $\{a_i,a_{i+3},b_i,b_{i+3},c_i,c_{i+3},d_i,d_{i+3}, e_i, e_{i+3}, f_i\}$ for $i\in\{0,1,2\}$;
    \item $\{s_i,t_i\}$ for $i \in \mathbb{Z}_6$; and
    \item $\{p_i,q_i\::\: i\in\Z_2\}\cup \{\infty\}$.
\end{itemize}
Blocks are obtained by developing modulo $6$ the subscripts of the noninfinite points in the base blocks given in Table~\ref{table.4gdd.2^6_5^1_11^3}.
The first three blocks in the first column generate two blocks each, while the last two blocks in that column generate three blocks each.
\end{proof}
\end{lemma}

\begin{table}[ht]
    \centering
    \caption{Base blocks for the $4$-GDD of type~$2^6 5^1 11^3$ in Lemma~\ref{lemma.4gdd.2^6_5^1_11^3}}
    \label{table.4gdd.2^6_5^1_11^3}
    \renewcommand{\arraystretch}{1.2}
    \[ \begin{array}{|l|l|l|l|l|}\hline
    \{c_{0},c_{2},c_{4},q_{0}\}  & \{a_{0},a_{1},c_{2},s_{0}\}  & \{a_{0},d_{4},f_{2},q_{0}\}  &  \{b_{0},c_{5},f_{1},p_{1}\} & \{c_{0},d_{2},f_{1},t_{2}\} \\\hline
    \{e_{0},e_{2},e_{4},q_{0}\}  & \{a_{0},a_{2},d_{1},t_{0}\}  & \{a_{0},e_{5},f_{1},s_{3}\}  &  \{b_{0},c_{2},t_{1},t_{2}\} & \{c_{0},q_{1},s_{2},t_{1}\} \\\hline
    \{s_{0},s_{2},s_{4},q_{0}\}  & \{a_{0},b_{1},c_{5},s_{4}\}  & \{a_{0},e_{2},q_{1},t_{2}\}  & \{b_{0},d_{5},e_{4},t_{0}\}  & \{d_{0},d_{1},p_{0},t_{3}\} \\\hline
    \{f_{0},f_{1},t_{0},t_{3}\}  & \{a_{0},b_{4},d_{2},s_{2}\}  & \{a_{0},e_{4},t_{1},t_{5}\}  & \{b_{0},e_{1},f_{2},s_{0}\}  & \{d_{0},e_{1},e_{2},s_{4}\}\  \\\hline
    \{f_{0},s_{0},s_{3},\infty\} & \{a_{0},b_{2},e_{1},p_{0}\}  & \{b_{0},b_{1},d_{2},q_{1}\}  &  \{c_{0},c_{1},e_{5},t_{4}\} & \{d_{0},s_{2},s_{3},t_{4}\} \  \\\hline
                                 & \{a_{0},b_{5},t_{3},\infty\} & \{b_{0},b_{2},s_{1},t_{5}\}  & \{c_{0},d_{1},d_{5},s_{0}\}  & \{e_{0},p_{0},s_{1},t_{4}\} \\\hline
                                 & \{a_{0},c_{4},p_{1},s_{1}\}  & \{b_{0},c_{1},e_{2},s_{2}\}  & \{c_{0},d_{4},e_{2},\infty\} & \\\hline
 \end{array} \]
\end{table}

\begin{lemma}
\label{lemma.4gdd.2^3_11^4}
There exists a $4$-GDD of type~$2^3 11^4$.

\begin{proof}
The points in the $4$-GDD are
\begin{itemize}
    \item $a_i,b_i,c_i,d_i,e_i,r_i,t_i$ for $i\in\Z_6$;
    \item $f_i$ for $i \in\Z_3$;
    \item $p_i,q_i$ for $i\in\Z_2$; and
    \item $\infty$.
\end{itemize}
The groups are
\begin{itemize}
    \item $\{a_i,a_{i+3},b_i,b_{i+3},c_i,c_{i+3},d_i,d_{i+3}, e_i, e_{i+3}, f_i\}$ and $\{t_i,t_{i+3}\}$ for $i\in\{0,1,2\}$; and
    \item $\{p_i,q_i\::\: i\in\Z_2\}\cup \{r_i\::\: i\in\Z_6\}\cup\{\infty\}$.
\end{itemize}
Blocks are obtained by developing modulo $6$ the subscripts of the noninfinite points in the base blocks given in Table~\ref{table.4gdd.2^3_11^4}.
The first two blocks in the first column generate two blocks each, while the third block in that column generates one block.
\end{proof}
\end{lemma}

\begin{table}[ht]
    \centering
    \caption{Base blocks for the $4$-GDD of type~$2^3 11^4$ in Lemma~\ref{lemma.4gdd.2^3_11^4}}
    \label{table.4gdd.2^3_11^4}
    \renewcommand{\arraystretch}{1.2}
    \[ \begin{array}{|l|l|l|l|l|l|}\hline
    \{a_{0},a_{2},a_{4},p_{0}\}  & \{a_{0},a_{1},c_{2},r_{0}\}  & \{a_{0},c_{4},f_{2},t_{2}\} & \{b_{0},b_{1},f_{2},r_{4}\} & \{c_{0},c_{1},r_{2},t_{0}\}  & \{d_{0},f_{2},r_{3},t_{4}\}\\\hline
    \{c_{0},c_{2},c_{4},q_{0}\}  & \{a_{0},b_{1},d_{2},p_{1}\}  & \{a_{0},d_{4},e_{2},r_{3}\} & \{b_{0},c_{4},e_{5},q_{1}\} & \{c_{0},d_{4},d_{5},t_{1}\}  & \{e_{0},e_{1},p_{0},t_{1}\}\\\hline
    \{f_{0},f_{1},f_{2},\infty\} & \{a_{0},b_{2},d_{1},q_{0}\}  & \{a_{0},d_{5},q_{1},t_{5}\} & \{b_{0},c_{2},e_{4},r_{1}\} & \{c_{0},d_{2},e_{4},\infty\} & \{e_{0},f_{1},q_{1},t_{2}\}\\\hline
                                 & \{a_{0},b_{5},e_{1},t_{0}\}  & \{a_{0},e_{5},f_{1},r_{4}\} & \{b_{0},c_{1},p_{1},t_{3}\} & \{c_{0},d_{1},f_{2},p_{1}\}  & \\\hline
                                 & \{a_{0},b_{4},t_{4},\infty\} & \{a_{0},r_{1},t_{1},t_{3}\} & \{b_{0},c_{5},r_{5},t_{2}\} & \{d_{0},d_{2},r_{2},t_{1}\}  & \\\hline
                                 & \{a_{0},c_{5},e_{4},r_{2}\}  & \{b_{0},b_{2},d_{4},r_{2}\} & \{b_{0},e_{1},t_{4},t_{5}\} & \{d_{0},e_{1},e_{5},r_{1}\}  & \\\hline
  \end{array} \]
\end{table}

\section{Acknowledgements} 
This research used the computational cluster Katana supported by Research Technology Services at UNSW Sydney.
The third author acknowledges the support from an Australian Government Research Training Program Scholarship and from the School of  Mathematics and Statistics, UNSW Sydney.
%The authors would like to thank the reviewers for detailed checking of the constructions and providing a number of useful comments.

\section*{ORCID}
R. J. R. Abel:     https://orcid.org/0000-0002-3632-9612\\
T. Britz:          https://orcid.org/0000-0003-4891-3055\\
Y. A. Bunjamin:    https://orcid.org/0000-0001-6849-2986\\
D. Combe:          https://orcid.org/0000-0002-1055-3894


\begin{thebibliography}{99}
%\frenchspacing

%\bibitem{mols1860}
%R.J.R. Abel,  Five MOLS of orders 18 and 60,  J. Combin. Des.  27 (2015), 135--139.

%\bibitem{ABCD}
%R.J.R.  Abel,   C.J. Colbourn and J.H. Dinitz, Mutually orthogonal Latin 
%squares (MOLS), in:  The CRC Handbook of Combinatorial Designs, Second Edition
%(C.  J.\ Colbourn and J.  H.\ Dinitz, eds.), CRC Press, Boca Raton FL, 2007, 160--193.

% \bibitem{ABBC.2t8s}
% R.J.R. Abel, T. Britz, Y.A. Bunjamin and  D. Combe, Group divisible designs with block size $4$ where the group sizes are congruent to $2 \pmod{3}$, Discrete Math. 345 (2022), 112740.

\bibitem{ABBC.2t5s}
R.J.R. Abel, T. Britz, Y.A. Bunjamin and  D. Combe, Group divisible designs with block size $4$ and group sizes $2$ and $5$,  J. Combin. Des. 30 (2022), 367--383.  

\bibitem{ABBCF.1mod3}
R.J.R. Abel, T. Britz, Y.A. Bunjamin, D. Combe and T. Feng, Group divisible designs with block size 4 where the group sizes are congruent to $1 \pmod{3}$, Discrete Math. 346 (2023), 113277.

\bibitem{ABC.30less}
R.J.R. Abel,  Y.A. Bunjamin and D. Combe, Some new group divisible designs with block size $4$ and two or three group sizes,  J. Combin. Des. 28 (2020), 614--628.  

\bibitem{ABC.50less}
R.J.R. Abel, Y.A. Bunjamin and  D. Combe, Existence of 4-GDDs with at most 50 points and $4$-GDDs of types $6^s 3^t$ and $9^s 3^t$, Discrete Math. 344 (2021), 112479.

% \bibitem{ABC.3562}
% R.J.R. Abel,  Y.A. Bunjamin and D. Combe, The 4-GDDs of type $3^5 6^2$, Discrete Math. 345 (2022), 112983.

%\bibitem{BR1979} A.E. Brouwer, Optimal packings of $K_4$'s into a %$K_n$. J. Combin. Theory Ser. A 26 (1979), no. 3, 278--297.


%\bibitem{BW}
%R.D. Baker  and R.M. Wilson, Nearly Kirkman triple systems,
%    Util. Math.  11 (1977), 289--296.

% \bibitem{BSH}
% A.E. Brouwer, A. Schrijver and H. Hanani,  Group divisible designs with block size 4, Discrete Math. 20 (1977), 1--10.


%\bibitem{dengreesshen3} D. Deng, R. Rees and H. Shen, Further results on nearly Kirkman triple systems with subsystems,  Discrete Math. 270, (2003),  99--114.

% \bibitem{Forbes2}
% A.D. Forbes,  Group divisible designs with block size 4 and type $g^u m^1$ II,  J. Combin. Des.  27 (2019), 311--349.

% \bibitem{Forbes3}
% A.D. Forbes,  Group divisible designs with block size 4 and type $g^u m^1$ III,  J. Combin. Des.  27 (2019), 643--672.

% \bibitem{Forbesgcc}
% A.D. Forbes,  Group divisible designs with block size 4 and type $g^u b^1 (gu/2)^1$,  Graphs Combin. 36 (2020), 1687--1703. % https://doi.org/10.1007/s00373-020-02213-5.

% \bibitem{Forbes1}
% A.D. Forbes and K.A. Forbes,  Group divisible designs with block size 4 and type $g^u m^1$,  J. Combin. Des.  26 (2018), 519--539.

% \bibitem{fmy}   
% S.C. Furino, Y. Miao and J. Yin, Frames and resolvable designs, CRC Press, Boca Raton FL, 1996.

% \bibitem{GeLing.2004.gum1}
% G. Ge and A. Ling, Group divisible designs with block size four and group type $g^u m^1$ for small $g$, Discrete Math. 285 (2004), 97--120.

% \bibitem{gerees}
% G. Ge and R.S. Rees, On group-divisible designs with block size four and group type $g^u m^1$, Des. Codes Cryptogr. 27 (2002), 5--24.

%\bibitem{ge0mod6}
% G. Ge and H. Wei,  Group divisible designs with block size 4  and group type 
% $g^u m^1$ for $g \equiv 0$ $($mod $6)$  J. Combin Des. 22  (2014), 26--52.  

% \bibitem{gezhu}
% G. Ge, R.S. Rees and L. Zhu, Group divisible designs with block size 4 and group type $g^u m^1$ with $m$ as small as possible,  J. Combin. Theory Ser. A  98 (2002), 357-376.  

% \bibitem{Hanani.1975}
% H. Hanani, Balanced incomplete block designs and related designs, Discrete Math. 11 (1975) 255–369.

%\bibitem{reesgrut}
%M. Gr\"{u}ttm\"{u}ller and R.S. Rees, Mandatory representation designs MRD$(4,k;v)$ with $k \equiv 1 \pmod{3}$, Util.  Math. 60 (2001), 153--180.   %R.S. Rees here

%   a\bibitem{Ingo}
% I. Janiszczak and R. Staszewski,  Isometry invariant permutation codes and mutually orthogonal
%  Latin squares,    J. Combin. Des.  27 (2019), 541--551.

% \bibitem{kreher.210note}
% D.L. Kreher, A.C.H. Ling, R.S. Rees and C.W.H. Lam, A note on $\{4\}$-GDDs of type $2^{10}$, Discrete Math. 261 (2003) 373–376.

\bibitem{krestin}
D. Kreher and D.R. Stinson, Small group divisible designs with block size four, J. Statist. Plann. Inference 58 (1997), 111--118.

% \bibitem{mckay.nauty}
% B.D. McKay and A. Piperno, Practical Graph Isomorphism, II,  J. Symbolic Comput. 60 (2014), 94--112.

% \bibitem{reesk=45}
% R.S. Rees, Group-divisible designs with block size $k$ having $k+1$ groups for $k=4,5$, J. Combin. Des. 8 (2000), 363--386. %R.S. Rees here

\bibitem{reesstinsubs}
R. Rees  and D.R. Stinson,  Kirkman triple systems with maximum subsystems, Ars Combin.  25 (1988), 125--132.

%\bibitem{reesstin89}   % at present unciteds; it contains 4-GDD$(3^{11} 9^1)$
%R. Rees  and D.R. Stinson,  On combinatorial designs with subdesigns, Discrete Math.  77 (1989), 259--279.

% \bibitem{reesstinonehole}
% R. Rees  and D.R. Stinson,  On the existence of incomplete designs of block size 4 having one hole, Util. Math. 35 (1989), 119--152.

%\bibitem{reesstinrgdd3}
%R. Rees and D.R. Stinson, On resolvable group-divisible designs with block size $3$, Ars Combin. 23 (1987), 107--120. % R. Rees not R.S. Rees in this title

% \bibitem{Schuster.gum1.mult8}
% E. Schuster, Group divisible designs with block size four and group type $g^u m^1$ where $g$ is a multiple of $8$, Discrete Math. 310 (2010), 2258--2270.

% \bibitem{schuster.2014.gum1}
% E. Schuster, New classes of group divisible designs with block size 4 and group type $g^u m^1$, J. Combin. Math. Combin. Comput. 91 (2014), 65--105.

%\bibitem{ICKPD4}
%L. Wang, R.J.R. Abel, D. Deng and J. Wang, Existence of incomplete canonical Kirkman packing designs, Discrete Math. 341 (2018), 536--554.

%\bibitem{wangshen}
%J. Wang and  H. Shen,  Existence of $(v, K_{1(3)} \cup \{w^*\})$-PBDs and its applications, Des. Codes Cryptogr. 46 (2008), 1--16.  

%\bibitem{geweikf}
%H. Wei and G. Ge,  Group divisible designs with block sizes from $K_{1(3)}$ and Kirkman frames of type $g^u m^1$,  Discrete Math.  329 (2014), 42--68.  

% \bibitem{WeiGe.2013.more.gum1}
% H. Wei and G. Ge, Group divisible designs with block size four and group type $g^u m^1$ for more small $g$, Discrete Math. 313 (2013), 2065--2083.

% \bibitem{WeiGe.4gdd.gpn1.g0mod6}
% H. Wei and G. Ge, Group divisible designs with block size four and group type $g^u m^1$ for $g \equiv 0 \pmod{6}$, J. Combin. Des. 22 (2014), 26--52.


% \bibitem{gegum}
% H. Wei and G. Ge, Group divisible designs with block size 4 and type $g^u m^1$, Des. Codes Cryptogr. 74 (2015), 243--282.

% \bibitem{Wilson}
% R.M. Wilson, Constructions and uses of pairwise balanced designs, Math. Centre Tracts  55 (1974), 18--41.

\end{thebibliography}
\end{document}